\documentclass{article}
\usepackage{amsfonts}
\usepackage{graphicx}
\usepackage{amsmath}
\usepackage[francais]{babel}

\newtheorem{theorem}{Theorem}

\newtheorem{proposition}[theorem]{Proposition}

\newtheorem{theoreme}[theorem]{Th\'eor\`eme}
\newtheorem{corollaire}[theorem]{Corollaire}

\newtheorem{remarque}[theorem]{Remarque}

\newenvironment{demonstration}[1][D\'emonstration]{\textbf{#1:} }
{\ \rule{0.5em}{0.5em}}
\input{tcilatex}

\begin{document}

\title{Champs de vecteurs et formes diff\'{e}rentielles sur une vari\'{e}t%
\'{e} des points proches}
\author{Basile Guy Richard BOSSOTO$^{(1)}$, Eug\`{e}ne OKASSA$^{(2)}$ \\
Universit\'{e} Marien NGOUABI, Facult\'{e} des Sciences,\\
D\'{e}partement de Math\'{e}matiques\\
B.P.69 - BRAZZAVILLE- (Congo)\\
e-mail: $^{(1)}$bossotob@yahoo.fr, $^{(2)}$eugeneokassa@yahoo.fr}
\date{}
\maketitle

\begin{abstract}
On consid\`{e}re $M$ une vari\'{e}t\'{e} diff\'{e}rentielle, $A$ une alg\`{e}%
bre locale au sens d'Andr\'{e} Weil, $M^{A}$ la vari\'{e}t\'{e} des points
proches de $M$ d'esp\`{e}ce $A$ et $\mathfrak{X}(M^{A})$ le module des
champs de vecteurs sur $M^{A}$. On donne une nouvelle d\'{e}finition des
champs de vecteurs sur $M^{A}$ et on montre que $\mathfrak{X}(M^{A})$ est
une alg\`{e}bre de Lie sur $A$. On \'{e}tudie la cohomologie des $A$-formes
diff\'{e}rentielles.

\textbf{Summary: }Let $M$ be a smooth manifold, $A$ a local algebra in sense
of Andr\'{e} Weil, $M^{A}$ the manifold of near points on $M$ of kind $A$
and $\mathfrak{X}(M^{A})$ the module of vector fields on $M^{A}$. We give a
new definition of vector fields on $M^{A}$ and we show that $\mathfrak{X}%
(M^{A})$ is a Lie algebra over $A$. We study the cohomology of $A$%
-differential forms.
\end{abstract}



\section{Introduction}

On consid\`{e}re une vari\'{e}t\'{e} lisse $M$, $A$ une alg\`{e}bre locale
(au sens d'Andr\'{e} Weil) et $M^{A}$ la vari\'{e}t\'{e} des points proches
de $M$ d'esp\`{e}ce $A$ $\cite{wei}$. Lorsque la vari\'{e}t\'{e} $M$ est de
dimension $n$, alors $M^{A}$ est une vari\'{e}t\'{e} lisse de dimension $%
n\cdot \dim (A)$.

On note $C^{\infty }(M)$ l'alg\`{e}bre des fonctions de classe $C^{\infty }$
sur $M$, $\mathfrak{X}(M)$ le $C^{\infty }(M)$-module des champs de vecteurs
sur $M$.

Lorsque $M$ et $N$ sont deux vari\'{e}t\'{e}s lisses et lorsque 
\begin{equation*}
h:M\longrightarrow N
\end{equation*}%
est une application diff\'{e}rentiable de classe $C^{\infty }$, alors
l'application 
\begin{equation*}
h^{A}:M^{A}\longrightarrow N^{A},\xi \longmapsto h^{A}(\xi ),
\end{equation*}%
telle que, pour tout $\varphi \in C^{\infty }(N)$, 
\begin{equation*}
\left[ h^{A}(\xi )\right] (\varphi )=\xi (\varphi \circ h)
\end{equation*}%
est diff\'{e}rentiable de classe $C^{\infty }$. Lorsque $h$ est un diff\'{e}%
omorphisme, il en est de m\^{e}me de $h^{A}$.

L'ensemble, $C^{\infty }(M^{A},A)$, des fonctions de classe $C^{\infty }$
sur $M^{A}$ \`{a} valeurs dans $A$, est une $A$-alg\`{e}bre commutative
unitaire.

En identifiant $%
\mathbb{R}
^{A}$ \`{a} $A$, pour $f\in C^{\infty }(M)$, l'application 
\begin{equation*}
f^{A}:M^{A}\longrightarrow A,\xi \longmapsto \xi (f)\text{,}
\end{equation*}%
est de classe $C^{\infty }$. De plus l'application%
\begin{equation*}
C^{\infty }(M)\longrightarrow C^{\infty }(M^{A},A),f\longmapsto f^{A},
\end{equation*}%
est un homomorphisme injectif d'alg\`{e}bres. Ainsi, on a:%
\begin{eqnarray*}
(f+g)^{A} &=&f^{A}+g^{A} \\
(\lambda \cdot f)^{A} &=&\lambda \cdot f^{A} \\
(f\cdot g)^{A} &=&f^{A}\cdot g^{A}
\end{eqnarray*}%
avec $\lambda \in \mathbb{R}$, $f$ et $g$ appartenant \`{a} $C^{\infty }(M)$.

Lorsque $(a_{\alpha })_{\alpha =1,2,...,\dim (A)}$ est une base de $A$ et
lorsque $(a_{\alpha }^{\ast })_{\alpha =1,2,...,\dim (A)}$ est la base duale
de la base $(a_{\alpha })_{\alpha =1,2,...,\dim (A)}$, l'application%
\begin{equation*}
\sigma :C^{\infty }(M^{A},A)\longrightarrow A\otimes C^{\infty
}(M^{A}),\varphi \longmapsto \dsum\limits_{\alpha =1}^{\dim (A)}a_{\alpha
}\otimes (a_{\alpha }^{\ast }\circ \varphi ),
\end{equation*}%
est un isomorphisme de $A$-alg\`{e}bres. Cet isomorphisme ne d\'{e}pend \'{e}%
videmment pas de la base choisie. L'application%
\begin{equation*}
\gamma :C^{\infty }(M)\longrightarrow A\otimes C^{\infty
}(M^{A}),f\longmapsto \sigma (f^{A}),
\end{equation*}%
est un morphisme d'alg\`{e}bres.

Dans toute la suite $M$ est une vari\'{e}t\'{e} lisse paracompacte de
dimension $n$.

Lorsque $(U,\varphi )$ est une carte locale de $M$ de syst\`{e}me de coordonn%
\'{e}es locales $(x_{1},x_{2},...,x_{n})$, l'application 
\begin{equation*}
U^{A}\longrightarrow A^{n},\xi \longmapsto (\xi (x_{1}),\xi (x_{2}),...,\xi
(x_{n})),
\end{equation*}%
est une bijection de $U^{A}$ sur un ouvert de $A^{n}$. On v\'{e}rifie que $%
M^{A}$ est une $A$-vari\'{e}t\'{e} de dimension $n$.

L'ensemble, $\mathfrak{X}(M^{A})$, des champs de vecteurs sur $M^{A}$ est 
\`{a} la fois un $C^{\infty }(M^{A})$-module et un $A$- module. Ce qui
signifie que $\mathfrak{X}(M^{A})$ est un $C^{\infty }(M^{A},A)$-module.

Dans ce travail, on \'{e}tudie la structure de $C^{\infty }(M^{A},A)$-module
de $\mathfrak{X}(M^{A})$. De cette nouvelle approche, on construit une
structure de $A$-alg\`{e}bre de Lie sur $\mathfrak{X}(M^{A})$, on d\'{e}%
finit les $A$-formes diff\'{e}rentielles et on en \'{e}tudie la cohomologie.

\section{Structure de $A$-alg\`{e}bre de Lie sur $\mathfrak{X}(M^{A})$}

\subsection{\protect\bigskip Vecteurs tangents sur $M^{A}$}

Pour $\xi \in M^{A}$, on note $T_{\xi }M^{A}$ l'espace tangent en $\xi \in
M^{A}$ et $Der_{\xi }\left[ C^{\infty }(M),A\right] $ l'ensemble des $\xi $-d%
\'{e}rivations de $C^{\infty }(M)$ dans $A$ c'est-\`{a}-dire l'ensemble des
applications $%
\mathbb{R}
$-lin\'{e}aires 
\begin{equation*}
v:C^{\infty }(M)\longrightarrow A
\end{equation*}%
telles que, pour $f$ et $g$ appartenant \`{a} $C^{\infty }(M)$, 
\begin{equation*}
v(fg)=v(f)\cdot \xi (g)+\xi (f)\cdot v(g)
\end{equation*}%
i.e.%
\begin{equation*}
v(fg)=v(f)\cdot g^{A}(\xi )+f^{A}(\xi )\cdot v(g)\text{.}
\end{equation*}

\begin{proposition}
\cite{oka1},\cite{oka2}L'application%
\begin{equation*}
T_{\xi }M^{A}=Der_{\xi }\left[ C^{\infty }(M^{A}),%
\mathbb{R}
\right] \longrightarrow Der_{\xi }\left[ C^{\infty }(M),A\right]
,v\longmapsto (id_{A}\otimes v)\circ \gamma ,
\end{equation*}%
est un isomorphisme d'espaces vectoriels .
\end{proposition}

Cet isomorphisme permet de transporter sur $T_{\xi }M^{A}$ la structure de $%
A $-module du $A$-module $Der_{\xi }\left[ C^{\infty }(M),A\right] $.

Ainsi:

\begin{corollaire}
Les assertions suivantes sont \'{e}quivalentes:

$1/$ $v$ est un vecteur tangent en $\xi \in M^{A}$;

$2/$ $v$ est une application $%
\mathbb{R}
$-lin\'{e}aire de $C^{\infty }(M)$ dans $A$ telle que, pour $f$ et $g$
appartenant \`{a} $C^{\infty }(M)$,%
\begin{equation*}
v(fg)=v(f)\cdot g^{A}(\xi )+f^{A}(\xi )\cdot v(g)\text{.}
\end{equation*}
\end{corollaire}

Lorsque $\xi \in M^{A}$, l'application 
\begin{equation*}
\widetilde{\xi }:C^{\infty }(M^{A},A)\longrightarrow A,\varphi \longmapsto
\varphi (\xi ),
\end{equation*}%
est un homomorphisme d'alg\`{e}bres. On note $Der_{\widetilde{\xi }}\left[
C^{\infty }(M^{A},A),A\right] $ le $A$-module des $\widetilde{\xi }$-d\'{e}%
rivations de $C^{\infty }(M^{A},A)$ dans $A$ c'est-\`{a}-dire l'ensemble des
applications $\mathbb{R}$-lin\'{e}aires 
\begin{equation*}
w:C^{\infty }(M^{A},A)\longrightarrow A
\end{equation*}%
telles que, pour $\varphi $ et $\psi $ appartenant \`{a} $C^{\infty
}(M^{A},A)$, 
\begin{equation*}
w(\varphi \cdot \psi )=w(\varphi )\cdot \widetilde{\xi }(\psi )+\widetilde{%
\xi }(\varphi )\cdot w(\psi )\text{.}
\end{equation*}

On d\'{e}duit le th\'{e}or\`{e}me suivant:

\begin{theoreme}
Si%
\begin{equation*}
v:C^{\infty }(M)\longrightarrow A
\end{equation*}%
est un vecteur tangent en $\xi \in M^{A}$, alors il existe une $\widetilde{%
\xi }$-d\'{e}rivation et une seule 
\begin{equation*}
\widetilde{v}:C^{\infty }(M^{A},A)\longrightarrow A
\end{equation*}%
telle que:

1/ $\widetilde{v}$ est $A$-lin\'{e}aire;

2/ $\ \widetilde{v}\left[ C^{\infty }(M^{A})\right] \subset \mathbb{R}$;

3/ $\ \widetilde{v}(f^{A})=v(f)$ pour tout $f\in C^{\infty }(M)$.
\end{theoreme}

\begin{demonstration}
Soit 
\begin{equation*}
v:C^{\infty }(M)\longrightarrow A
\end{equation*}
un vecteur tangent en $\xi \in M^{A}$ et soit 
\begin{equation*}
\overline{v}:C^{\infty }(M^{A})\longrightarrow 
\mathbb{R}%
\end{equation*}%
l'unique d\'{e}rivation telle que 
\begin{equation*}
(id_{A}\otimes \overline{v})\circ \gamma =v\text{.}
\end{equation*}%
L'application%
\begin{equation*}
\widetilde{v}=(id_{A}\otimes \overline{v})\circ \sigma :C^{\infty
}(M^{A},A)\longrightarrow A
\end{equation*}%
r\'{e}pond \`{a} la question.
\end{demonstration}

\subsection{Champs de vecteurs sur $M^{A}$}

On note $Der_{\gamma }\left[ C^{\infty }(M),A\otimes C^{\infty }(M^{A})%
\right] $ le $A\otimes C^{\infty }(M^{A})$-module des $\gamma $-d\'{e}%
riva-tions de $C^{\infty }(M)$ dans $A\otimes C^{\infty }(M^{A})$ i.e.
l'ensemble des applications $%
\mathbb{R}
$-lin\'{e}aires 
\begin{equation*}
\varphi :C^{\infty }(M)\longrightarrow A\otimes C^{\infty }(M^{A})
\end{equation*}%
telles que, pour $f$ et $g$ appartenant \`{a} $C^{\infty }(M)$,\qquad\ 
\begin{equation*}
\varphi (fg)=\varphi (f)\cdot \gamma (g)+\gamma (f)\cdot \varphi (g)\text{.}
\end{equation*}

Une d\'{e}rivation de $C^{\infty }(M)$ dans $C^{\infty }(M^{A},A)$ est une
application $%
\mathbb{R}
$- lin\'{e}aire 
\begin{equation*}
Y:C^{\infty }(M)\longrightarrow C^{\infty }(M^{A},A)
\end{equation*}%
telle que, pour $f$ et $g$ appartenant \`{a} $C^{\infty }(M)$,%
\begin{equation*}
Y(fg)=Y(f)\cdot g^{A}+f^{A}\cdot Y(g)\text{.}
\end{equation*}%
Ainsi une d\'{e}rivation de $C^{\infty }(M)$ dans $C^{\infty }(M^{A},A)$ est
une d\'{e}rivation par rapport \`{a} l'homomorphisme 
\begin{equation*}
C^{\infty }(M)\longrightarrow C^{\infty }(M^{A},A),f\longmapsto f^{A}\text{.}
\end{equation*}%
Il s'ensuit que l'ensemble, $Der\left[ C^{\infty }(M),C^{\infty }(M^{A},A)%
\right] $, des d\'{e}rivations de $C^{\infty }(M)$ dans $C^{\infty
}(M^{A},A) $ est un $C^{\infty }(M^{A},A)$-module.

\begin{proposition}
\cite{oka1},\cite{oka2}L'application 
\begin{equation*}
Der\left[ C^{\infty }(M^{A})\right] \longrightarrow Der_{\gamma }\left[
C^{\infty }(M),A\otimes C^{\infty }(M^{A})\right] ,X\longmapsto
(id_{A}\otimes X)\circ \gamma ,
\end{equation*}%
est un isomorphisme de $C^{\infty }(M^{A})$-modules.
\end{proposition}

Il s'ensuit:

\begin{corollaire}
L'application%
\begin{equation*}
Der\left[ C^{\infty }(M^{A})\right] \longrightarrow Der\left[ C^{\infty
}(M),C^{\infty }(M^{A},A)\right] ,X\longmapsto \sigma ^{-1}\circ
(id_{A}\otimes X)\circ \gamma ,
\end{equation*}%
est un isomorphisme de $C^{\infty }(M^{A})$-modules.
\end{corollaire}

Cet isomorphisme permet de transporter sur $Der\left[ C^{\infty }(M^{A})%
\right] $ la structure de $C^{\infty }(M^{A},A)$-module de $Der\left[
C^{\infty }(M),C^{\infty }(M^{A},A)\right] $.

Ainsi:

\begin{corollaire}
Les assertions suivantes sont \'{e}quivalentes:

$1/$ Un champ de vecteurs sur $M^{A}$ est une section diff\'{e}rentiable du
fibr\'{e} tangent $(TM^{A},\pi _{M^{A}},M^{A})$;

$2/$ Un champ de vecteurs sur $M^{A}$ est une d\'{e}rivation de $C^{\infty
}(M^{A})$;

$3/$ Un champ de vecteurs sur $M^{A}$ est une d\'{e}rivation de $C^{\infty
}(M)$ dans $C^{\infty }(M^{A},A)$.
\end{corollaire}

On d\'{e}duit le th\'{e}or\`{e}me suivant:

\begin{theoreme}
Si $X$ est un champ de vecteurs sur $M^{A}$ consid\'{e}r\'{e} comme d\'{e}%
rivation de $C^{\infty }(M)$ dans $C^{\infty }(M^{A},A)$, alors il existe
une d\'{e}rivation et une seule 
\begin{equation*}
\widetilde{X}:C^{\infty }(M^{A},A)\longrightarrow C^{\infty }(M^{A},A)
\end{equation*}%
telle que

$1/$ $\widetilde{X}$ est $A$-lin\'{e}aire;

$2/$ $\widetilde{X}\left[ C^{\infty }(M^{A})\right] \subset C^{\infty
}(M^{A})$;

$3/\ \widetilde{X}(f^{A})=X(f)$ pour tout $f\in C^{\infty }(M)$.
\end{theoreme}

\begin{demonstration}
Si $X$ est un champ de vecteurs sur $M^{A}$ consid\'{e}r\'{e} comme d\'{e}%
rivation de $C^{\infty }(M)$ dans $C^{\infty }(M^{A},A)$ et si 
\begin{equation*}
\overline{X}:C^{\infty }(M^{A})\longrightarrow C^{\infty }(M^{A})
\end{equation*}%
est l'unique d\'{e}rivation telle que%
\begin{equation*}
\sigma ^{-1}\circ (id_{A}\otimes \overline{X})\circ \gamma =X\text{,}
\end{equation*}%
alors l'application%
\begin{equation*}
\widetilde{X}=\sigma ^{-1}\circ (id_{A}\otimes \overline{X})\circ \sigma
:C^{\infty }(M^{A},A)\longrightarrow C^{\infty }(M^{A},A)
\end{equation*}%
r\'{e}pond \`{a} la question.
\end{demonstration}

\begin{remarque}
Si $X$ est un champ de vecteurs sur $M^{A}$ consid\'{e}r\'{e} comme d\'{e}%
rivation de $C^{\infty }(M)$ dans $C^{\infty }(M^{A},A)$, alors $\widetilde{X%
}$ s'annule sur $A$.
\end{remarque}

\begin{proposition}
Si $\mu :A\longrightarrow A$ est un endomorphisme, $f\in C^{\infty }(M)$ et $%
X:C^{\infty }(M)$ $\longrightarrow $ $C^{\infty }(M^{A},A)$ un champ de
vecteurs sur $M^{A}$, alors 
\begin{equation*}
\widetilde{X}(\mu \circ f^{A})=\mu \circ X(f)\text{.}
\end{equation*}
\end{proposition}

\begin{demonstration}
De $\widetilde{X}(f^{A})=X(f)$, on a 
\begin{equation*}
\widetilde{X}\left[ \dsum\limits_{\alpha =1}^{\dim (A)}(a_{\alpha }^{\ast
}\circ f^{A})\cdot a_{\alpha }\right] =\dsum\limits_{\alpha =1}^{\dim
(A)}(a_{\alpha }^{\ast }\circ X(f)\cdot a_{\alpha }\text{.}
\end{equation*}

Ce qui donne 
\begin{equation*}
\dsum\limits_{\alpha =1}^{\dim (A)}\widetilde{X}(a_{\alpha }^{\ast }\circ
f^{A})\cdot a_{\alpha }=\dsum\limits_{\alpha =1}^{\dim (A)}(a_{\alpha
}^{\ast }\circ X(f)\cdot a_{\alpha }\text{.}
\end{equation*}

Ainsi $\widetilde{X}(a_{\alpha }^{\ast }\circ f^{A})=(a_{\alpha }^{\ast
}\circ X(f)$ pour tout $(a_{\alpha }^{\ast })_{i=1,2,...,\dim (A)}$.Comme 
\begin{equation*}
\mu \circ f^{A}=\dsum\limits_{\alpha =1}^{\dim (A)}(a_{\alpha }^{\ast }\circ
f^{A})\cdot \mu (a_{\alpha })\text{,}
\end{equation*}

on d\'{e}duit que%
\begin{eqnarray*}
\widetilde{X}(\mu \circ f^{A}) &=&\dsum\limits_{\alpha =1}^{\dim (A)}%
\widetilde{X}(a_{\alpha }^{\ast }\circ f^{A})\cdot \mu (a_{\alpha }) \\
&=&\dsum\limits_{\alpha =1}^{\dim (A)}(a_{\alpha }^{\ast }\circ X(f)\cdot
\mu (a_{\alpha }) \\
&=&\mu \circ X(f)\text{.}
\end{eqnarray*}

D'o\`{u} l'assertion.
\end{demonstration}

\begin{theoreme}
Si $X$ et $Y$ sont deux champs de vecteurs sur $M$\bigskip $^{A}$ consid\'{e}%
r\'{e}s comme d\'{e}rivations de $C^{\infty }(M)$ dans $C^{\infty }(M^{A},A)$%
, alors le crochet 
\begin{equation*}
\left[ X,Y\right] =\widetilde{X}\circ Y-\widetilde{Y}\circ X:C^{\infty
}(M)\longrightarrow C^{\infty }(M^{A},A)
\end{equation*}%
est un champ de vecteurs sur $M^{A}$.
\end{theoreme}

\begin{demonstration}
L'application est manifestement $%
\mathbb{R}
$-lin\'{e}aire. Pour $f$ et $g$ appartenant \`{a} $C^{\infty }(M)$, on a 
\begin{eqnarray*}
\left[ X,Y\right] (fg) &=&\widetilde{X}\left[ Y(fg)\right] -\widetilde{Y}%
\left[ X(fg)\right] \\
&=&\widetilde{X}\left[ Y(f)\cdot g^{A}+f^{A}\cdot Y(g))\right] \\
&&-\widetilde{Y}\left[ X(f)\cdot g^{A}+f^{A}\cdot X(g)\right] \\
&=&\widetilde{X}\left[ Y(f)\right] \cdot g^{A}+Y(f)\cdot \widetilde{X}%
(g^{A})+\widetilde{X}(f^{A})\cdot Y(g)+f^{A}\cdot \widetilde{X}\left[ Y(g)%
\right] \\
&&-\widetilde{Y}\left[ X(f)\right] \cdot g^{A}-X(f)\cdot \widetilde{Y}%
(g^{A})-\widetilde{Y}(f^{A})\cdot X(g)-f^{A}\cdot \widetilde{Y}\left[ X(g)%
\right] \\
&=&\widetilde{X}\left[ Y(f)\right] \cdot g^{A}+Y(f)\cdot X(g)+X(f)\cdot
Y(g)+f^{A}\cdot \widetilde{X}\left[ Y(g)\right] \\
&&-\widetilde{Y}\left[ X(f)\right] \cdot g^{A}-X(f)\cdot Y(g)-Y(f)\cdot
X(g)-f^{A}\cdot \widetilde{Y}\left[ X(g)\right] \\
&=&(\widetilde{X}\left[ Y(f)\right] -\widetilde{Y}\left[ X(f)\right] )\cdot
g^{A}+f^{A}\cdot (\widetilde{X}\left[ Y(g)\right] -\widetilde{Y}\left[ X(g)%
\right] ) \\
&=&(\widetilde{X}\circ Y-\widetilde{Y}\circ X)(f)\cdot g^{A}+f^{A}\cdot (%
\widetilde{X}\circ Y-\widetilde{Y}\circ X)(g) \\
&=&\left[ X,Y\right] (f)\cdot g^{A}+f^{A}\cdot \left[ X,Y\right] (g)\text{.}
\end{eqnarray*}

D'o\`{u} l'assertion.
\end{demonstration}

\begin{proposition}
Si $X$ et $Y$ sont deux champs de vecteurs sur $M$\bigskip $^{A}$ consid\'{e}%
r\'{e}s comme d\'{e}rivations de $C^{\infty }(M)$ dans $C^{\infty }(M^{A},A)$
et si $\varphi \in C^{\infty }(M^{A},A)$, alors%
\begin{equation*}
\left[ \widetilde{X},\widetilde{Y}\right] =\widetilde{\left[ X,Y\right] }
\end{equation*}%
et 
\begin{equation*}
\widetilde{\varphi \cdot X}=\varphi \cdot \widetilde{X}\text{.}
\end{equation*}
\end{proposition}

\begin{demonstration}
Pour $f\in C^{\infty }(M)$, on a\textbf{\ }%
\begin{eqnarray*}
\left[ \widetilde{X},\widetilde{Y}\right] (f^{A}) &=&\widetilde{X}\left[ 
\widetilde{Y}(f^{A})\right] -\widetilde{Y}\left[ \widetilde{X}(f^{A})\right]
\\
&=&\widetilde{X}\left[ Y(f)\right] -\widetilde{Y}\left[ X(f)\right] \\
&=&(\widetilde{X}\circ Y-\widetilde{Y}\circ X)(f) \\
&=&\left[ X,Y\right] (f)\text{.}
\end{eqnarray*}%
Comme $\widetilde{\left[ X,Y\right] }$ est l'unique d\'{e}rivation de $%
C^{\infty }(M^{A},A)$ telle que $\widetilde{\left[ X,Y\right] }(f^{A})=\left[
X,Y\right] (f)$ pour tout $f\in C^{\infty }(M)$, on d\'{e}duit que 
\begin{equation*}
\left[ \widetilde{X},\widetilde{Y}\right] =\widetilde{\left[ X,Y\right] }%
\text{.}
\end{equation*}

De m\^{e}me 
\begin{eqnarray*}
(\varphi \cdot \widetilde{X})(f^{A}) &=&\varphi \cdot (\widetilde{X})(f^{A})
\\
&=&\varphi \cdot X(f) \\
&=&(\varphi \cdot X)(f)\text{.}
\end{eqnarray*}

Comme $\widetilde{\varphi \cdot X}$ est l'unique d\'{e}rivation de $%
C^{\infty }(M^{A},A)$ telle que $($ $\widetilde{\varphi \cdot X}%
)(f^{A})=(\varphi \cdot X)(f)$ pour tout $f\in C^{\infty }(M)$, on d\'{e}%
duit que 
\begin{equation*}
\widetilde{\varphi \cdot X}=\varphi \cdot \widetilde{X}\text{.}
\end{equation*}

D'o\`{u} les deux assertions.
\end{demonstration}

\begin{proposition}
Si $\varphi \in C^{\infty }(M^{A},A)$, si $X$ et $Y$ sont deux champs de
vecteurs sur $M$\bigskip $^{A}$ consid\'{e}r\'{e}s comme d\'{e}rivations de $%
C^{\infty }(M)$ dans $C^{\infty }(M^{A},A)$, alors 
\begin{equation*}
\left[ X,\varphi \cdot Y\right] =\widetilde{X}(\varphi )\cdot Y+\varphi
\cdot \left[ X,Y\right] \text{.}
\end{equation*}
\end{proposition}

La d\'{e}monstration ne pr\'{e}sente aucune difficult\'{e}.

\begin{theoreme}
L'application 
\begin{equation*}
\mathfrak{X}(M^{A})\times \mathfrak{X}(M^{A})\longrightarrow \mathfrak{X}%
(M^{A}),(X,Y)\longmapsto \left[ X,Y\right] ,
\end{equation*}%
est $A$-bilin\'{e}aire altern\'{e}e et d\'{e}finit une structure de $A$-alg%
\`{e}bre de Lie sur $\mathfrak{X}(M^{A})$.
\end{theoreme}

\begin{demonstration}
Lorsque $X$ est un champ de vecteurs sur $M$\bigskip $^{A}$ consid\'{e}r\'{e}
comme d\'{e}rivation de $C^{\infty }(M)$ dans $C^{\infty }(M^{A},A)$ et
lorsque $a\in A$, on a%
\begin{equation*}
\left[ X,a\cdot Y\right] =\widetilde{X}(a)\cdot Y+a\cdot \left[ X,Y\right] 
\text{.}
\end{equation*}%
Comme $\widetilde{X}$ s'annule sur $A$, il s'ensuit que l'application 
\begin{equation*}
\mathfrak{X}(M^{A})\times \mathfrak{X}(M^{A})\longrightarrow \mathfrak{X}%
(M^{A}),(X,Y)\longmapsto \left[ X,Y\right]
\end{equation*}

est $A$-bilin\'{e}aire altern\'{e}e.

Pour tous champs de vecteurs $X$, $Y$, $Z$ sur $M^{A}$ consid\'{e}r\'{e}s
comme d\'{e}rivations de $C^{\infty }(M)$ dans $C^{\infty }(M^{A},A)$, on a:

\begin{eqnarray*}
&&\left[ X,\left[ Y,Z\right] \right] +\left[ Y,\left[ Z,X\right] \right] +%
\left[ Z,\left[ X,Y\right] \right] \\
&=&\widetilde{X}\circ \left[ Y,Z\right] -\widetilde{\left[ Y,Z\right] }\circ
X \\
&&+\widetilde{Y}\circ \left[ Z,X\right] -\widetilde{\left[ Z,X\right] }\circ
Y \\
&&+\widetilde{Z}\circ \left[ X,Y\right] -\widetilde{\left[ X,Y\right] }\circ
Z \\
&=&\widetilde{X}\circ (\widetilde{Y}\circ Z-\widetilde{Z}\circ Y)-\left[ 
\widetilde{Y},\widetilde{Z}\right] \circ X \\
&&+\widetilde{Y}\circ (\widetilde{Z}\circ X-\widetilde{X}\circ Z)-\left[ 
\widetilde{Z},\widetilde{X}\right] \circ Y \\
&&+\widetilde{Z}\circ (\widetilde{X}\circ Y-\widetilde{Y}\circ X)-\left[ 
\widetilde{X},\widetilde{Y}\right] \circ Z \\
&=&\widetilde{X}\circ \widetilde{Y}\circ Z-\widetilde{X}\circ \widetilde{Z}%
\circ Y-\widetilde{Y}\circ \widetilde{Z}\circ X+\widetilde{Z}\circ 
\widetilde{Y}\circ X \\
&&+\widetilde{Y}\circ \widetilde{Z}\circ X-\widetilde{Y}\circ \widetilde{X}%
\circ Z-\widetilde{Z}\circ \widetilde{X}\circ Y+\widetilde{X}\circ 
\widetilde{Z}\circ Y \\
&&+\widetilde{Z}\circ \widetilde{X}\circ Y-\widetilde{Z}\circ \widetilde{Y}%
\circ X-\widetilde{X}\circ \widetilde{Y}\circ Z+\widetilde{Y}\circ 
\widetilde{X}\circ Z \\
&=&0\text{.}
\end{eqnarray*}

D'o\`{u} l'assertion.
\end{demonstration}

\begin{remarque}
En consid\'{e}rant $\mathfrak{X}(M^{A})$ uniquement comme module sur $%
C^{\infty }(M^{A})$, $\mathfrak{X}(M^{A})$ ne peut \^{e}tre une alg\`{e}bre
de Lie sur $A$.
\end{remarque}

\begin{corollaire}
L'application 
\begin{equation*}
\mathfrak{X}(M^{A})\longrightarrow Der\left[ C^{\infty }(M^{A},A)\right]
,X\longmapsto \widetilde{X},
\end{equation*}%
est \`{a} la fois un morphisme de $C^{\infty }(M^{A},A)$-modules et un
morphisme de $A$-alg\`{e}bres de Lie.
\end{corollaire}

\subsubsection{\textbf{Prolongements \`{a} }$M^{A}$\textbf{\ des champs de
vecteurs sur }$M$}

\begin{proposition}
Si%
\begin{equation*}
\theta :C^{\infty }(M)\longrightarrow C^{\infty }(M)
\end{equation*}%
est un champ de vecteurs sur $M$, alors l'application 
\begin{equation*}
\theta ^{A}:C^{\infty }(M)\longrightarrow C^{\infty }(M^{A},A),f\longmapsto 
\left[ \theta (f)\right] ^{A},
\end{equation*}%
est un champ de vecteurs sur $M^{A}$.
\end{proposition}

\begin{demonstration}
L'application $\theta ^{A}$ est manifestement $%
\mathbb{R}
$-lin\'{e}aire. Pour $f$ et $g$ appartenant \`{a} $C^{\infty }(M)$, on a:%
\begin{eqnarray*}
\theta ^{A}(fg) &=&\left[ \theta (fg)\right] ^{A} \\
&=&\left[ \theta (f)\cdot g+f\cdot \theta (g)\right] ^{A} \\
&=&\left[ \theta (f)\right] ^{A}\cdot g^{A}+f^{A}\cdot \left[ \theta (g)%
\right] ^{A} \\
&=&\theta ^{A}(f)\cdot g^{A}+f^{A}\cdot \theta ^{A}(g)\text{.}
\end{eqnarray*}

Ainsi $\theta ^{A}$ est un champ de vecteurs sur $M^{A}$.
\end{demonstration}

On dit que le champ de vecteurs $\theta ^{A}$ est le prolongement \`{a} $%
M^{A}$ du champ de vecteurs $\theta $ sur $M$.

\begin{proposition}
Si $\theta $, $\theta _{1}$ et $\theta _{2}$ sont des champs de vecteurs sur 
$M$ et si $f\in C^{\infty }(M)$, alors \textbf{\ }%
\begin{eqnarray*}
(\theta _{1}+\theta _{2})^{A} &=&\theta _{1}^{A}+\theta _{2}^{A}\text{;} \\
(f\cdot \theta )^{A} &=&f^{A}\cdot \theta ^{A}\text{;} \\
\widetilde{(f\cdot \theta )^{A}} &=&f^{A}\cdot \widetilde{\theta ^{A}}\text{;%
} \\
\left[ \theta _{1}^{A},\theta _{2}^{A}\right] &=&\left[ \theta _{1},\theta
_{2}\right] ^{A}\text{.}
\end{eqnarray*}%
et l'application 
\begin{equation*}
\mathfrak{X}(M)\longrightarrow Der\left[ C^{\infty }(M^{A},A)\right] ,\theta
\longmapsto \widetilde{\theta ^{A}},
\end{equation*}%
est un homomorphisme d'alg\`{e}bres de Lie r\'{e}elles.
\end{proposition}

La d\'{e}monstration ne pr\'{e}sente aucune difficult\'{e}.

\subsubsection{\textbf{Champs de vecteurs sur }$M^{A}$\textbf{\ provenant
des d\'{e}rivations de }$A$}

\begin{proposition}
Si $d$ est une d\'{e}rivation de $A$, alors l'application \ 
\begin{equation*}
d^{\ast }:C^{\infty }(M)\longrightarrow C^{\infty }(M^{A},A),f\longmapsto
(-d)\circ f^{A},
\end{equation*}%
est un champ de vecteurs sur $M^{A}$.
\end{proposition}

\begin{demonstration}
On v\'{e}rifie que l'application $d^{\ast }$ est $%
\mathbb{R}
$-lin\'{e}aire. Pour $f$ et $g$ appartenant \`{a} $C^{\infty }(M)$ et pour $%
\xi \in M^{A}$, on a:%
\begin{eqnarray*}
d^{\ast }(fg)(\xi ) &=&(-d)\circ (fg)^{A}(\xi ) \\
&=&(-d)\circ (f^{A}\cdot g^{A})(\xi ) \\
&=&(-d)\left[ f^{A}(\xi )\cdot g^{A}(\xi )\right] \  \\
&=&(-d)\left[ f^{A}(\xi )\right] \cdot g^{A}(\xi )+f^{A}(\xi )\cdot (-d)%
\left[ g^{A}(\xi )\right] \\
&=&\left[ (-d)\circ f^{A}\right] (\xi )\cdot g^{A}(\xi )+f^{A}(\xi )\cdot %
\left[ (-d)\circ f^{A}\right] (\xi ) \\
&=&(\left[ (-d)\circ f^{A}\right] \cdot g^{A}+f^{A}\cdot \left[ (-d)\circ
f^{A}\right] )(\xi ) \\
&=&\left[ d^{\ast }(f)\cdot g^{A}+f^{A}\cdot d^{\ast }(g)\right] (\xi )\text{%
.}
\end{eqnarray*}

Comme $\xi $ est quelconque, on d\'{e}duit que 
\begin{equation*}
d^{\ast }(fg)=d^{\ast }(f)\cdot g^{A}+f^{A}\cdot d^{\ast }(g)\text{.}
\end{equation*}

Ainsi, $d^{\ast }$ est un champ de vecteurs sur $M^{A}$.
\end{demonstration}

On dit que le champ de vecteurs $d^{\ast }$ est le champ de vecteurs sur $%
M^{A}$ associ\'{e} \`{a} la d\'{e}rivation $d$ de $A$.

On a les r\'{e}sultats suivants:

\begin{proposition}
Si $d_{1}$, $d_{2}$ , $d$ sont trois d\'{e}rivations de $A$, $a$ un \'{e}l%
\'{e}ment de $A$ et $\theta :C^{\infty }(M)\longrightarrow C^{\infty }(M)$
un champ de vecteurs sur $M$, alors%
\begin{eqnarray*}
\left[ d_{1}^{\ast },d_{2}^{\ast }\right] &=&\left[ d_{1},d_{2}\right]
^{\ast }\text{;} \\
(a\cdot d)^{\ast } &=&a\cdot d^{\ast }\text{;} \\
\left[ d^{\ast },\theta ^{A}\right] &=&0\text{.}
\end{eqnarray*}
\end{proposition}

\begin{demonstration}
$\bigskip $La d\'{e}monstration des deux premi\`{e}res assertions ne pr\'{e}%
sente aucune difficult\'{e}.

Pour la derni\`{e}re assertion, lorsque $f$ $\in $ $C^{\infty }(M)$ on a%
\begin{eqnarray*}
\left[ d^{\ast },\theta ^{A}\right] (f) &=&(\widetilde{d^{\ast }}\circ
\theta ^{A}-\widetilde{\theta ^{A}}\circ d^{\ast })(f) \\
&=&(\widetilde{d^{\ast }}\circ \theta ^{A})(f)-(\widetilde{\theta ^{A}}\circ
d^{\ast })(f) \\
&=&(\widetilde{d^{\ast }})\left[ \theta ^{A}(f)\right] -(\widetilde{\theta
^{A}})\left[ d^{\ast }(f)\right] \\
&=&(\widetilde{d^{\ast }})(\left[ \theta (f)\right] ^{A})-(\widetilde{\theta
^{A}})\left[ d^{\ast }(f)\right] \\
&=&d^{\ast }\left[ \theta (f)\right] +(\widetilde{\theta ^{A}})\left[ d\circ
f^{A}\right] \\
&=&(-d)\circ \left[ \theta (f)\right] ^{A}+(\widetilde{\theta ^{A}})\left[
d\circ f^{A}\right] \text{.}
\end{eqnarray*}

Compte tenu de la proposition $9$, on a%
\begin{equation*}
(\widetilde{\theta ^{A}})\left[ d\circ f^{A}\right] =d\circ \theta ^{A}(f)%
\text{.}
\end{equation*}

Ainsi 
\begin{eqnarray*}
\left[ d^{\ast },\theta ^{A}\right] (f) &=&(-d)\circ \left[ \theta (f)\right]
^{A}+(\widetilde{\theta ^{A}})\left[ d\circ f^{A}\right] \\
&=&(-d)\circ \theta ^{A}(f)+d\circ \theta ^{A}(f) \\
&=&0\text{.}
\end{eqnarray*}

Comme $f$ est quelconque, on d\'{e}duit que $\left[ d^{\ast },\theta ^{A}%
\right] =0$.
\end{demonstration}

\section{$A$-formes diff\'{e}rentielles}

Un $A$-covecteur en $\xi \in M^{A}$ est une forme lin\'{e}aire sur le $A$%
-module $T_{\xi }M^{A}$. L'ensemble, $T_{\xi }^{\ast }M^{A}$, des $A$%
-covecteurs en $\xi \in M^{A}$ est un $A$-module libre de dimension $n$ et 
\begin{equation*}
T^{\ast }M^{A}=\bigcup\limits_{\xi \in M^{A}}T_{\xi }^{\ast }M^{A}
\end{equation*}%
est une $A$-vari\'{e}t\'{e} de dimension $2n$. L'ensemble, $\Lambda
^{1}(M^{A},A)$, des sections diff\'{e}rentiables de $T^{\ast }M^{A}$ est un $%
C^{\infty }(M^{A},A)$-module et on dit que $\Lambda ^{1}(M^{A},A)$ est le $%
C^{\infty }(M^{A},A)$-module des $A$-formes diff\'{e}rentielles de degr\'{e} 
$1$.

Pour $p\in \mathbb{N}$ et pour $\xi \in M^{A}$, on note $\mathcal{L}%
_{alt}^{p}(T_{\xi }M^{A},A)$ le $A$-module des formes multilin\'{e}aires
altern\'{e}es de degr\'{e} $p$ sur le $A$-module $T_{\xi }M^{A}$. On a \'{e}%
videmment 
\begin{equation*}
\mathcal{L}_{alt}^{0}(T_{\xi }M^{A},A)=A\text{.}
\end{equation*}%
Comme dans le cas r\'{e}el, pour deux entiers $p$ et $q$, on d\'{e}finit le
produit ext\'{e}rieur 
\begin{equation*}
\Lambda :\mathcal{L}_{alt}^{p}(T_{\xi }M^{A},A)\times \mathcal{L}%
_{alt}^{q}(T_{\xi }M^{A},A)\longrightarrow \mathcal{L}_{alt}^{p+q}(T_{\xi
}M^{A},A),(\alpha ,\beta )\longmapsto \alpha \Lambda \beta .
\end{equation*}%
L'ensemble 
\begin{equation*}
A^{p}(T^{\ast }M^{A},A)=\bigcup\limits_{\xi \in M^{A}}\mathcal{L}%
_{alt}^{p}(T_{\xi }M^{A},A)
\end{equation*}%
est une $A$-vari\'{e}t\'{e} de dimension $n+C_{n}^{p}$. L'ensemble, $\Lambda
^{p}(M^{A},A)$, des sections diff\'{e}rentiables de $A^{p}(T^{\ast }M^{A},A)$
est un $C^{\infty }(M^{A},A)$-module. On dit que $\Lambda ^{p}(M^{A},A)$ est
le $C^{\infty }(M^{A},A)$-module des $A$-formes diff\'{e}rentielles de degr%
\'{e} $p$ sur $M^{A}$ et que 
\begin{equation*}
\Lambda (M^{A},A)=\bigoplus\limits_{p=0}^{n}\Lambda ^{p}(M^{A},A)
\end{equation*}%
est l'alg\`{e}bre des $A$-formes diff\'{e}rentielles sur $M^{A}$. L'alg\`{e}%
bre $\Lambda (M^{A},A)$ des $A$-formes diff\'{e}rentielles sur $M^{A}$ est
canoniquement isomorphe \`{a} $A\otimes \Lambda (M^{A})$. On a 
\begin{equation*}
\Lambda ^{0}(M^{A},A)=C^{\infty }(M^{A},A)\text{.}
\end{equation*}

\begin{theoreme}
\cite{mor},\cite{oka3}Si $\eta $ est une $A$-forme diff\'{e}rentielle de degr%
\'{e} $p$ sur $M^{A}$, alors il existe une $A$-forme diff\'{e}rentielle de
degr\'{e} $p$ et une seule 
\begin{equation*}
\eta ^{A}:\mathfrak{X}(M^{A})\times \mathfrak{X}(M^{A})\times ...\times 
\mathfrak{X}(M)^{A})\longrightarrow C^{\infty }(M^{A},A)
\end{equation*}%
telle que , pour $p$ champs de vecteurs $\theta _{1},\theta _{2},...,\theta
_{p}$ sur $M$ et pour $p$ fonctions $f_{1},f_{2},...,f_{p}$ sur $M$, 
\begin{equation*}
\eta ^{A}(f_{1}^{A}\cdot \theta _{1}^{A},f_{2}^{A}\cdot \theta
_{2}^{A},...,f_{p}^{A}\theta _{p}^{A})=f_{1}^{A}\cdot f_{2}^{A}\cdot
...\cdot f_{p}^{A}\cdot \left[ \eta (\theta _{1},\theta _{2},...,\theta _{p})%
\right] ^{A}\text{.}
\end{equation*}
\end{theoreme}

Lorsque $\eta $ est une forme diff\'{e}rentielle sur $M$, la $A$-forme diff%
\'{e}rentielle $\eta ^{A}$ est le prolongement \`{a} $M^{A}$ de la forme diff%
\'{e}rentielle $\eta $.

\subsection{La $d^{A}$- cohomologie}

L'application 
\begin{equation*}
\Lambda (M)\longrightarrow \Lambda (M^{A},A),\omega \longmapsto \omega ^{A},
\end{equation*}%
est un morphisme de $%
\mathbb{R}
$-alg\`{e}bres gradu\'{e}es.

Si 
\begin{equation*}
d:\Lambda (M)\longrightarrow \Lambda (M)
\end{equation*}%
est l'op\'{e}rateur de diff\'{e}rentiation ext\'{e}rieure, en notant 
\begin{equation*}
d^{A}:\Lambda (M^{A},A)\longrightarrow \Lambda (M^{A},A)
\end{equation*}%
l'op\'{e}rateur de cohomologie associ\'{e} \`{a} la repr\'{e}sentation 
\begin{equation*}
\mathfrak{X}(M^{A})\longrightarrow Der\left[ C^{\infty }(M^{A},A)\right]
,X\longmapsto \widetilde{X}\text{.}
\end{equation*}

\begin{proposition}
L'application 
\begin{equation*}
d^{A}:\Lambda (M^{A},A)\longrightarrow \Lambda (M^{A},A)
\end{equation*}%
est $A$-lin\'{e}aire et 
\begin{equation*}
d^{A}(\omega ^{A})=(d\omega )^{A}
\end{equation*}%
pour tout $\omega \in \Lambda (M)$.
\end{proposition}

\begin{demonstration}
On v\'{e}rifie que $d^{A}$ est $A$-lin\'{e}aire. Si $\omega \in \Lambda
^{p}(M)$, pour $\theta _{1},\theta _{2},...,\theta _{p+1}$ champs de
vecteurs sur $M$, on a%
\begin{eqnarray*}
&&\left[ d^{A}(\omega ^{A})\right] (\theta _{1}^{A},\theta
_{2}^{A},...,\theta _{p+1}^{A}) \\
&=&\overset{p+1}{\underset{i=1}{\sum }}(-1)^{i-1}\widetilde{\theta _{i}^{A}}%
\left[ \omega ^{A}(\theta _{1}^{A},\theta _{2}^{A},...,\widehat{\theta
_{i}^{A}},...,\theta _{p+1}^{A})\right] \\
&&+\underset{1\leq i<j\leq p+1}{\sum }(-1)^{i+j}\ \omega ^{A}(\left[ \theta
_{i}^{A},\theta _{j}^{A}\right] ,\theta _{1}^{A},...,\widehat{\theta _{i}^{A}%
}\ ...,\widehat{\theta _{j}^{A}},...,\theta _{p+1}^{A}) \\
&=&\overset{p+1}{\underset{i=1}{\sum }}(-1)^{i-1}\widetilde{\theta _{i}^{A}}%
\left[ (\omega (\theta _{1},\theta _{2},...,\widehat{\theta _{i}},...,\theta
_{p+1}))^{A}\right] \\
&&+\underset{1\leq i<j\leq p+1}{\sum }(-1)^{i+j}\ (\omega (\left[ \theta
_{i},\theta _{j}\right] ,\theta _{1},...,\widehat{\theta _{i}}\ ...,\widehat{%
\theta _{j}},...,\theta _{p+1}))^{A} \\
&=&\overset{p+1}{\underset{i=1}{\sum }}(-1)^{i-1}\theta _{i}^{A}\left[
\omega (\theta _{1},\theta _{2},...,\widehat{\theta _{i}},...,\theta _{p+1})%
\right] \\
&&+\underset{1\leq i<j\leq p+1}{\sum }(-1)^{i+j}\ (\omega (\left[ \theta
_{i},\theta _{j}\right] ,\theta _{1},...,\widehat{\theta _{i}}\ ...,\widehat{%
\theta _{j}},...,\theta _{p+1}))^{A} \\
&=&\overset{p+1}{\underset{i=1}{\sum }}(-1)^{i-1}\left( \theta _{i}\left[
\omega (\theta _{1},\theta _{2},...,\widehat{\theta _{i}},...,\theta _{p+1})%
\right] \right) ^{A} \\
&&+\underset{1\leq i<j\leq p+1}{\sum }(-1)^{i+j}\ (\omega (\left[ \theta
_{i},\theta _{j}\right] ,\theta _{1},...,\widehat{\theta _{i}}\ ...,\widehat{%
\theta _{j}},...,\theta _{p+1}))^{A} \\
&=&(d\omega )^{A}(\theta _{1}^{A},\theta _{2}^{A},...,\theta _{p+1}^{A}) \\
&=&\left[ (d\omega )(\theta _{1},\theta _{2},...,\theta _{p+1})\right] ^{A}%
\text{.}
\end{eqnarray*}

Compte tenu du th\'{e}or\`{e}me $20$, on d\'{e}duit que $d^{A}(\omega
^{A})=(d\omega )^{A}$.
\end{demonstration}

L'application 
\begin{equation*}
A\times \Lambda (M)\longrightarrow \Lambda (M^{A},A),(a,\omega )\longmapsto
a\cdot \omega ^{A}
\end{equation*}%
est $%
\mathbb{R}
$-bilin\'{e}aire et induit un morphisme du complexe diff\'{e}rentiel $%
(A\otimes \Lambda (M),id_{A}\otimes d)$ dans le complexe diff\'{e}rentiel $%
(\Lambda (M^{A},A),d^{A})$.

On note $H_{dR}(M)$ la cohomologie de de Rham de la vari\'{e}t\'{e} diff\'{e}%
rentielle $M$ et $H(M^{A},A)$ la cohomologie du complexe diff\'{e}rentiel $%
(\Lambda (M^{A},A),d^{A})$.

On dit que $H(M^{A},A)$ est la $d^{A}$-cohomologie sur la vari\'{e}t\'{e}
des points proches $M^{A}$. Les espaces $A\otimes H_{dR}^{p}(M^{A})$ et $%
H^{p}(M^{A},A)$ sont canoniquement isomorphes.

En particulier si la vari\'{e}t\'{e} $M^{A}$ est connexe, alors l'espace $%
H^{0}(M^{A},A)$ s'identifie canoniquement \`{a} $A$.


\begin{thebibliography}{9}
\bibitem{mor} A. Morimoto, Prolongation of connections to bundles of
infinitely near points, J.Diff.Geom., t.11, 1976, p. 479-498.

\bibitem{oka1} E. Okassa, Prolongements des champs de vecteurs \`{a} des vari%
\'{e}t\'{e}s des points proches, C.R. Acad. Sc. Paris, t.300, S\'{e}rie I, n$%
{{}^\circ}%
$ 6, 1985, p.173-176.

\bibitem{oka2} E. Okassa, Prolongement des champs de vecteurs \`{a} des vari%
\'{e}t\'{e}s des points proches, Annales Facult\'{e} des Sciences de
Toulouse, Vol. VIII, n$%
{{}^\circ}%
$ 3, 1986-1987, 349-366.

\bibitem{oka3} E. Okassa, Rel\`{e}vements des structures symplectiques et
pseudo-riemanniennes \`{a} des vari\'{e}t\'{e}s des points proches, Nagoya
Math. J., Vol.115 (1989), 63-71.

\bibitem{yan1} K. Yano and S. Ishihara, Tangent and cotangent bundles, Diff.
Geom. Marcel Dekker, New-York, 1973.

\bibitem{yan2} K. Yano and E.M. Patterson, Vertical and complete lifts from
a manifold to its cotangent bundles, Jour. Math. Soc. Japan, 19 (1967),
91-113.

\bibitem{wei} A. Weil, Th\'{e}orie des points proches sur les vari\'{e}t\'{e}%
s diff\'{e}rentiables, Colloque G\'{e}om. Diff. Strasbourg, 1953, 111-117.
\end{thebibliography}
\end{document}